# ON TATE-SHAFAREVICH GROUPS OF ABELIAN VARIETIES

Cristian D. Gonzalez-Avilés

ABSTRACT.
   Let $K/F$ be a finite Galois extension of number fields with Galois group $G$, let $A$ be an abelian variety defined over $F$, and let $\text{III}(A_{/K})$ and $\text{III}(A_{/F})$ denote, respectively, the Tate-Shafarevich groups of $A$ over $K$ and of $A$ over $F$. Assuming that these groups are finite, we derive, under certain restrictions on $A$ and $K/F$, a formula for the order of the subgroup of $\text{III}(A_{/K})$ of $G$-invariant elements. As a corollary, we obtain a simple formula relating the orders of $\text{III}(A_{/K})$, $\text{III}(A_{/F})$ and $\text{III}(A^{\chi}_{/F})$ when $K/F$ is a quadratic extension and $A^\chi$ is the twist of $A$ by the non-trivial character $\chi$ of $G$.

## 1. Introduction

This paper is the first progress report of an ongoing investigation whose aim is to determine the behavior of the Tate-Shafarevich group of an abelian variety $A$ under extensions of the field of definition of $A$. To be precise, let $A$ be an abelian variety defined over a number field $F$, let $K/F$ be a finite Galois extension with Galois group $G$, and let $\text{III}(A_{/K})$ and $\text{III}(A_{/F})$ denote, respectively, the Tate-Shafarevich groups of $A$ over $K$ and of $A$ over $F$. We assume throughout that these groups are finite. Then our chief aim is to find a simple relation between the orders of $\text{III}(A_{/K})$ and $\text{III}(A_{/F})$, if such a relation exists. A partial solution to this problem is implicit in a 1972 paper of Milne ([9], Corollary to Theorem 3), who obtained his result making certain assumptions on $\text{End}_K(A) \otimes_{\mathbb{Z}} \mathbb{Q}$. We have adopted a different approach here, which works well for abelian varieties $A$ and field extensions $K/F$ as above which satisfy the following two conditions:

(A) $\hat{H}^p(G, A(K)) = \hat{H}^p(G, A'(K)) = 0$ for all $p$.
(B) Either $F$ is totally imaginary or both $A(\mathbb{R})$ and $A'(\mathbb{R})$ are connected.

Here $A'$ denotes the dual abelian variety of $A$. Thus for example $A$ could be an elliptic curve over $\mathbb{Q}$ given by a Weierstrass equation of negative discriminant and $K$ could be a finite Galois extension of $\mathbb{Q}$ such that $A(K)$ is finite and of order prime to $\#G = [K : \mathbb{Q}]$ (see Corollary V.2.3.1 of [14] and §6 of [1]). Our main result is the following.

**Main Theorem.** *Assume that conditions (A) and (B) above hold. Then*

$$\#\text{III}(A_{/K})^G = \#\text{III}(A_{/F}) \cdot \prod_{v \in S} \#H^1(G_w, A(K_w)),$$

1991 *Mathematics Subject Classification.* 11G40, 11G05.
Supported by Fondecyt grant 1981175.

Typeset by $\mathcal{A}_{\mathcal{M}}\mathcal{S}$-TEX

where: $S$ denotes the set of primes of $F$ obtained by collecting together the primes that ramify in $K/F$ and the primes of bad reduction for $A_{/F}$, $w$ is a fixed prime of $K$ lying above $v$ for each $v \in S$, and $G_w$ denotes $\mathrm{Gal}(K_w/F_v)$. In addition,

$$\#H^1(G, \mathrm{III}(A_{/K})) = \prod_{v \in S} \#H^2(G_w, A(K_w)).$$

The above theorem has the following corollary, which solves the problem of relating $\#\mathrm{III}(A_{/K})$ to $\#\mathrm{III}(A_{/F})$ in a special case.

**Corollary.** *Suppose that $K/F$ is a quadratic extension and let $\chi$ denote the non-trivial character of $G = \mathrm{Gal}(K/F)$. Assume that conditions (A) and (B) above hold for both $A$ and its quadratic twist $A^\chi$. Then*

$$\#\mathrm{III}(A_{/K}) = \#\mathrm{III}(A_{/F}) \cdot \#\mathrm{III}(A^\chi_{/F}) \cdot \prod_{v \in S} \#H^1(G_w, A(K_w)).$$

## Acknowledgements

It is a pleasure to acknowledge the help rendered me by Jean-Louis Colliot-Thélène, who kindly provided the proof of Proposition 4.2 below. I also thank James S. Milne for some helpful remarks and David Rohrlich for his encouragement while I wrote this paper.

## 2. Local Computations

If $M$ is a topological abelian group, we will write $M^*$ for the group of continuous characters of finite order of $M$, i.e. $M^* = \mathrm{Hom}_{\mathrm{cts}}(M, \mathbb{Q}/\mathbb{Z})$. Also, if $G$ is a finite group, $M$ is a $G$-module and $p$ is any integer, $\hat{H}^p(G, M)$ will denote the $p$-th Tate cohomology group of $M$ (see §6 of [1]). In particular if we write $M_G$ for the largest quotient of $M$ on which $G$ acts trivially and $N^\star : M_G \to M^G$ for the map induced by multiplication by $N = \sum_{\sigma \in G} \sigma \in \mathbb{Z}[G]$ on $M$, then

(1) $$\hat{H}^{-1}(G, M) = \ker(N^\star) \quad \text{and} \quad \hat{H}^0(G, M) = \mathrm{coker}(N^\star).$$

Let $A$ be an abelian variety defined over a number field $F$. For any field $L \supset F$, we will write $G_L$ for $\mathrm{Gal}(\bar{L}/L)$, where $\bar{L}$ is an algebraic closure of $L$. Further, we will write $H^p(L, A)$ for $H^p(G_L, A(\bar{L}))$ and $A'$ for the dual abelian variety of $A$.

Now let $K$ be a finite Galois extension of $F$ and let $G$ be the Galois group of $K$ over $F$. For any prime $w$ of $K$, we let $G_w = \mathrm{Gal}(K_w/F_v) \subset G$ be the decomposition group of $w$ over $F$, where $v$ is the prime of $F$ lying below $w$. Finally, we will write $\mathrm{res}_w$ for the local restriction map $H^1(F_v, A) \to H^1(K_w, A)$.

**Lemma 2.1.** *Let $w$ be a prime of $K$ and let $v$ be the prime of $F$ lying below $w$.*
*(i) If $w$ is archimedean, then there is an exact sequence*

$$0 \to H^1(G_w, A(K_w)) \to H^1(F_v, A) \xrightarrow{\mathrm{res}_w} H^1(K_w, A)^{G_w} \to 0.$$

(ii) If $w$ is non-archimedean, then there is an exact sequence

$$0 \to H^1(G_w, A(K_w)) \to H^1(F_v, A) \xrightarrow{\text{res}_w} H^1(K_w, A)^{G_w} \to H^2(G_w, A(K_w)) \to 0.$$

*Proof.* Assertion (i) is easy to check. Assertion (ii) follows from the exactness of the sequence

$$H^1(G_w, A(K_w)) \hookrightarrow H^1(F_v, A) \xrightarrow{\text{res}_w} H^1(K_w, A)^{G_w} \to H^2(G_w, A(K_w)) \to H^2(F_v, A)$$

(which is the exact sequence of terms of low degree belonging to the Hochschild-Serre spectral sequence $H^p(G_w, H^q(K_w, A)) \Longrightarrow H^{p+q}(F_v, A)$) and the fact that $H^2(F_v, A) = 0$ for $v$ non-archimedean (see [5] and Corollary I.3.4 of [8]). □

In what follows, $H^0(F_v, A')$ denotes $A'(F_v)$ unless $v$ is archimedean, in which case it denotes $\hat{H}^0(F_v, A') = A'(F_v)/N_{\bar{F}_v/F_v}A'(\bar{F}_v)$. Similarly for $H^0(K_w, A')$.

**Lemma 2.2.** *Let $v$ be any prime of $F$. Then the dual of the map $\bigoplus_{w|v} \text{res}_w : H^1(F_v, A) \to \bigoplus_{w|v} H^1(K_w, A)$ is the map $\prod_{w|v} H^0(K_w, A') \to H^0(F_v, A')$ induced by*

$$\prod_{w|v} A'(K_w) \to A'(F_v), \quad (x_w)_{w|v} \mapsto \sum_{w|v} N_{K_w/F_v}(x_w).$$

*Proof.* If $v$ is archimedean, the verification of the above statement is straightforward, using Remark I.3.7 of [8]. If $v$ is non-archimedean, the lemma follows easily from Tate's local duality theory [15]. □

Now let $S$ denote the set of primes of $F$ obtained by collecting together all primes which ramify in $K/F$ and all primes of bad reduction for $A_{/F}$. Further, let $S_\infty$ be the set of archimedean primes of $F$.

**Lemma 2.3.** *Let $w$ be a prime of $K$ and let $v$ be the prime of $F$ lying below $w$. Assume that $v \notin S \cup S_\infty$. Then for every $p \geq 1$,*

$$H^p(G_w, A(K_w)) = 0.$$

*Proof.* The case $p = 1$ of this result is well-known ([7], Corollary 4.4). For the general case, see Lemma 3.5 of [12]. □

Recall $G = \text{Gal}(K/F)$ and let $v$ be any prime of $F$. Then $\bigoplus_{w|v} H^q(K_w, A)$ can be made into a $G$-module in the following natural way. For $\sigma \in G$ and $(\xi_w)_{w|v} \in \bigoplus_{w|v} H^q(K_w, A)$, let $\sigma(\xi_w)_{w|v} = (\sigma_*^{-1}\xi_{\sigma w})_{w|v}$, where $\sigma_* : H^q(K_w, A) \to H^q(K_{\sigma w}, A)$ is the homomorphism associated to the maps $G_{K_{\sigma w}} \to G_{K_w}$, $\nu \mapsto \bar{\sigma}^{-1}\nu\bar{\sigma}$, and $A(\bar{K}_w) \to A(\bar{K}_{\sigma w})$, $P \mapsto \bar{\sigma}P$, where $\bar{\sigma} : \bar{K}_w \xrightarrow{\sim} \bar{K}_{\sigma w}$ is some lifting of $\sigma : K_w \xrightarrow{\sim} K_{\sigma w}$ (see [13], p. 115). It is not difficult to see that with this $G$-action, $\bigoplus_{w|v} H^q(K_w, A)$ becomes a semi-local $G$-module in the sense of [3]. Thus we have the following

**Lemma 2.4.** *Let $v$ be any prime of $F$. Then for every $p \geq 0$ and $q \geq 0$, there is a canonical isomorphism*

$$H^p\left(G, \bigoplus_{w|v} H^q(K_w, A)\right) \simeq H^p(G_w, H^q(K_w, A)),$$

*where the $w$ on the right denotes any prime of $K$ lying above $v$.*

*Proof.* See §2.1 of [3]. □

Let $v$ be a prime of $F$. It is easy to check that the image of the map $\bigoplus_{w|v} \mathrm{res}_w : H^1(F_v, A) \to \bigoplus_{w|v} H^1(K_w, A)$ is actually contained in $(\bigoplus_{w|v} H^1(K_w, A))^G$. Thus we have a map

$$\mathrm{res} : \bigoplus_v H^1(F_v, A) \to \left(\bigoplus_w H^1(K_w, A)\right)^G = \bigoplus_v \left(\bigoplus_{w|v} H^1(K_w, A)\right)^G,$$

namely $\mathrm{res} = \bigoplus_v \bigoplus_{w|v} \mathrm{res}_w$. Now recall the sets $S$ and $S_\infty$ defined above.

**Proposition 2.5.** *There are canonical isomorphisms*

$$\ker(\mathrm{res}) \simeq \bigoplus_{v \in S \cup S_\infty} H^1(G_w, A(K_w))$$

$$\mathrm{coker}(\mathrm{res}) \simeq \bigoplus_{v \in S} H^2(G_w, A(K_w)),$$

*where $w$ denotes a fixed prime of $K$ lying above $v$ for each $v \in S \cup S_\infty$.*

*Proof.* It suffices to compute, for any $v$, the kernel and cokernel of $s \circ \bigoplus_{w|v} \mathrm{res}_w$, where $s : (\bigoplus_{w|v} H^1(K_w, A))^G \to H^1(K_w, A)^{G_w}$ is the semi-local isomorphism of Lemma 2.4 corresponding to $p = 0$ and $q = 1$. Now the effect of $s$ is simply to project onto the $w$ coordinate (see [3]), from which it follows that $s \circ \bigoplus_{w|v} \mathrm{res}_w = \mathrm{res}_w$. The proposition now follows from Lemmas 2.1 and 2.3. □

## 3. Global Computations

Recall $G = \mathrm{Gal}(K/F)$. We will write $H^2(G, A(K))_{\mathrm{tr}}$ for the kernel of the natural inflation map $H^2(G, A(K)) \to H^2(F, A)$.

**Lemma 3.1.** *Let $\mathrm{Res} : H^1(F, A) \to H^1(K, A)^G$ be the global restriction map. Then*

$$\ker(\mathrm{Res}) \simeq H^1(G, A(K)) \quad \text{and} \quad \mathrm{coker}(\mathrm{Res}) \simeq H^2(G, A(K))_{\mathrm{tr}}.$$

*Proof.* This follows from the exactness of the sequence

$$0 \to H^1(G, A(K)) \to H^1(F, A) \xrightarrow{\mathrm{Res}} H^1(K, A)^G \to H^2(G, A(K)) \xrightarrow{\mathrm{inf}} H^2(F, A),$$

which is the exact sequence of terms of low degree belonging to the Hochschild-Serre spectral sequence $H^p(G, H^q(K, A)) \implies H^{p+q}(F, A)$. See [5]. □

In the next lemma, we write $r_1(F)$ for the number of real archimedean primes of $F$ and $A(\mathbb{R})^\circ$ (resp. $A'(\mathbb{R})^\circ$) for the identity component of $A(\mathbb{R})$ (resp. $A'(\mathbb{R})$).

**Lemma 3.2.** *Suppose that $q \geq 2$. If $q$ is even (resp. odd) then $H^q(F, A)$ is isomorphic to the direct sum of $r_1(F)$ copies of $A(\mathbb{R})/A(\mathbb{R})^\circ$ (resp. $A'(\mathbb{R})/A'(\mathbb{R})^\circ$).*

Proof. By Theorem I.6.26(c) of [8], the localization homomorphism $H^q(F, A) \to \bigoplus_{v \text{ real}} H^q(F_v, A)$ is an isomorphism. On the other hand, Remark I.3.7 of [8] shows that $H^q(\mathbb{R}, A)$ is isomorphic to either $\hat{H}^0(\mathbb{R}, A) = A(\mathbb{R})/A(\mathbb{R})^\circ$ if $q$ is even or to $\hat{H}^0(\mathbb{R}, A') = A'(\mathbb{R})/A'(\mathbb{R})^\circ$ if $q$ is odd. The lemma is now immediate. □

Let $\text{Ш}(A_{/K})$ and $\text{Ш}(A_{/F})$ denote the Tate-Shafarevich groups of $A$ over $K$ and of $A$ over $F$, respectively. These groups are defined by the exactness of the sequences

$$0 \to \text{Ш}(A_{/F}) \to H^1(F, A) \xrightarrow{\lambda_F} \bigoplus_v H^1(F_v, A) \to \text{coker}(\lambda_F) \to 0$$

and

$$0 \to \text{Ш}(A_{/K}) \to H^1(K, A) \xrightarrow{\lambda_K} \bigoplus_w H^1(K_w, A) \to \text{coker}(\lambda_K) \to 0,$$

where $\lambda_F$ and $\lambda_K$ are the natural localization maps. In what follows, we will assume true the well-known conjecture that $\text{Ш}(A_{/K})$ and $\text{Ш}(A_{/F})$ are finite groups. This conjecture has been verified in some special cases by Rubin [11] and Kolyvagin [6].

In the statement of the next proposition, we view $A'(K)$ and $A'(F)$ as topological groups with the profinite topology.

**Proposition 3.3.** *There are canonical $G$-isomorphisms*

$$\text{coker}(\lambda_K) \simeq A'(K)^* \quad \text{and} \quad \text{coker}(\lambda_F) \simeq A'(F)^*.$$

*Proof.* This follows from the finiteness of $\text{Ш}(A_{/K})$ and $\text{Ш}(A_{/F})$. The isomorphism $\text{coker}(\lambda_K) \simeq A'(K)^*$ is induced by the map $\bigoplus_w H^1(K_w, A) \to A'(K)^*$ which is dual to the diagonal embedding $A'(K) \to \prod_w H^0(K_w, A')$, and similarly for $\text{coker}(\lambda_F)$. See Theorem I.6.13 and Remark I.6.14 of [8]. □

Recall the map $\text{res}: \bigoplus_v H^1(F_v, A) \to \left(\bigoplus_w H^1(K_w, A)\right)^G$ introduced in §2. We have the following commutative diagram

$$\begin{array}{ccc} H^1(F, A) & \xrightarrow{\lambda_F} & \bigoplus_v H^1(F_v, A) \\ {\scriptstyle \text{Res}}\downarrow & & {\scriptstyle \text{res}}\downarrow \\ H^1(K, A)^G & \xrightarrow{\lambda_K} & \left(\bigoplus_w H^1(K_w, A)\right)^G. \end{array}$$

It follows that the map res induces a map

$$\rho : \text{coker}(\lambda_F) \to \text{coker}(\lambda_K)^G.$$

**Proposition 3.4.** *There are canonical isomorphisms*

$$\ker(\rho) \simeq \hat{H}^0(G, A'(K))^* \quad \text{and} \quad \text{coker}(\rho) \simeq \hat{H}^{-1}(G, A'(K))^*.$$

Proof. Let $N_{K/F} : A'(K) \to A'(K)$ be the global norm map. Then for any prime $v$ of $F$, $N_{K/F} = \sum_{w|v} N_{K_w/F_v}$, where $N_{K_w/F_v}$ denotes, for each $w|v$, the local norm map $A'(K_w) \to A'(K_w)$ (see Theorem I.15.3 of [10]). Thus we have a commutative diagram

$$\begin{array}{ccc} \prod_v H^0(F_v, A') & \longleftarrow & A'(F) \\ \uparrow & & \uparrow \\ \prod_v \left(\prod_{w|v} H^0(K_w, A')\right)_G & \longleftarrow & A'(K)_G, \end{array}$$

in which the horizontal maps are induced by the diagonal embeddings $A'(F) \to \prod_v H^0(F_v, A')$ and $A'(K) \to \prod_w H^0(K_w, A')$, the $v$-component of the left-hand vertical map is induced by the map $\prod_{w|v} H^0(K_w, A') \to H^0(F_v, A')$ of Lemma 2.2, and the right-hand vertical map is the map $N^\star_{K/F} : A'(K)_G \to A'(F) = A'(K)^G$ induced by $N_{K/F}$. The dual of the above diagram is the commutative diagram

$$\begin{array}{ccc} \bigoplus_v H^1(F_v, A) & \longrightarrow & A'(F)^* \\ \downarrow & & \downarrow \\ \left(\bigoplus_w H^1(K_w, A)\right)^G & \longrightarrow & (A'(K)^*)^G, \end{array}$$

where, by Lemma 2.2, the left-hand vertical map is the map res. It follows that under the isomorphisms $\mathrm{coker}(\lambda_K) \simeq A'(K)^*$ and $\mathrm{coker}(\lambda_F) \simeq A'(F)^*$ of Proposition 3.3, the map $\rho : \mathrm{coker}(\lambda_F) \to \mathrm{coker}(\lambda_K)^G$ corresponds to the dual of $N^\star_{K/F}$ (see the proof of Proposition 3.3). The lemma now follows easily from formula (1) of §2. □

## 4. The main result

We now make the following two assumptions on the abelian variety $A$ and field extension $K/F$ we are considering. These assumptions will remain in force for the rest of the paper.

(A) $\hat{H}^p(G, A(K)) = \hat{H}^p(G, A'(K)) = 0$ for all $p$.
(B) Either $F$ is totally imaginary or both $A(\mathbb{R})$ and $A'(\mathbb{R})$ are connected.

**Lemma 4.1.**

(i) For every archimedean prime $w$ of $K$, $H^1(G_w, A(K_w)) = 0$.
(ii) For all $q \geq 2$, $H^q(F, A) = H^q(K, A) = 0$.

*Proof.* Both assertions follow from assumption (B) above. See the statement and proof of Lemma 3.2. □

**Proposition 4.2.** *For every $p \geq 1$,*

$$H^p(G, H^1(K, A)) = 0.$$

*Proof.* By Lemma 4.1(ii), $H^q(K, A) = 0$ for all $q \geq 2$. Therefore by Theorem XV.5.11 of [2], there is an infinite exact sequence

$$\cdots \to H^{p+1}(F, A) \to H^p(G, H^1(K, A)) \to H^{p+2}(G, A(K)) \to H^{p+2}(F, A) \to \cdots.$$

The proposition now follows from Lemma 4.1(ii) and assumption (A) above. □

Now consider the commutative diagram with exact rows

$$\begin{array}{ccccccccc} 0 & \longrightarrow & \mathrm{im}(\lambda_F) & \longrightarrow & \bigoplus_v H^1(F_v, A) & \longrightarrow & \mathrm{coker}(\lambda_F) & \longrightarrow & 0 \\ & & \mathrm{res}' \downarrow & & \mathrm{res} \downarrow & & \rho \downarrow & & \\ 0 & \longrightarrow & \mathrm{im}(\lambda_K)^G & \longrightarrow & \left(\bigoplus_w H^1(K_w, A)\right)^G & \longrightarrow & \mathrm{coker}(\lambda_K)^G, & & \end{array}$$

where the maps res and $\rho$ are as defined previously, and res' is induced by res. Applying the snake lemma to this diagram yields the exact sequence

$$0 \to \ker(\text{res}') \to \ker(\text{res}) \to \ker(\rho) \to \text{coker}(\text{res}') \to \text{coker}(\text{res}) \to \text{coker}(\rho).$$

Now since $\ker(\rho) = \text{coker}(\rho) = 0$ by Proposition 3.4 and assumption (A), we conclude that there are isomorphisms

$$\ker(\text{res}') \simeq \ker(\text{res}) \quad \text{and} \quad \text{coker}(\text{res}') \simeq \text{coker}(\text{res}).$$

**Proposition 4.3.** *There are canonical isomorphisms*

$$\ker(\text{res}') \simeq \bigoplus_{v \in S} H^1(G_w, A(K_w))$$

$$\text{coker}(\text{res}') \simeq \bigoplus_{v \in S} H^2(G_w, A(K_w)),$$

*where $w$ denotes a fixed prime of $K$ lying above $v$ for each $v \in S$.*

*Proof.* This follows from the preceding discussion and Proposition 2.5 together with Lemma 4.1(i). □

**Theorem 4.4.** *Assuming conditions (A) and (B) above, we have*

$$\#\Sha(A_{/K})^G = \#\Sha(A_{/F}) \cdot \prod_{v \in S} \#H^1(G_w, A(K_w)),$$

*where $S$ denotes the set of primes of $F$ obtained by collecting together the primes that ramify in $K/F$ and the primes of bad reduction for $A_{/F}$. Furthermore,*

$$\#H^1(G, \Sha(A_{/K})) = \prod_{v \in S} \#H^2(G_w, A(K_w)).$$

*Proof.* Consider the commutative diagram with exact rows

$$\begin{array}{ccccccccc} 0 & \longrightarrow & \Sha(A_{/F}) & \longrightarrow & H^1(F, A) & \longrightarrow & \text{im}(\lambda_F) & \longrightarrow & 0 \\ & & \text{Res}' \downarrow & & \text{Res} \downarrow & & \text{res}' \downarrow & & \\ 0 & \longrightarrow & \Sha(A_{/K})^G & \longrightarrow & H^1(K, A)^G & \longrightarrow & \text{im}(\lambda_K)^G & \longrightarrow & \ldots, \end{array}$$

in which the bottom row is the long $G$-cohomology sequence associated with the exact sequence $0 \to \Sha(A_{/K}) \to H^1(K, A) \to \text{im}(\lambda_K) \to 0$, Res' is induced by the global restriction map Res, and res' is as defined above. Applying the snake lemma to the above diagram yields the exact sequence

$$0 \to \ker(\text{Res}') \to \ker(\text{Res}) \to \ker(\text{res}') \to \text{coker}(\text{Res}') \to \text{coker}(\text{Res}) \to$$
$$\text{coker}(\text{res}') \to H^1(G, \Sha(A_{/K})) \to H^1(G, H^1(K, A)).$$

Now Lemma 3.1 together with assumption (A) yields $\ker(\text{Res}) = \text{coker}(\text{Res}) = 0$, while $H^1(G, H^1(K, A)) = 0$ by Proposition 4.2. It follows that Res' is injective with cokernel isomorphic to the kernel of res', and that $H^1(G, \Sha(A_{/K})) \simeq \text{coker}(\text{res}')$. The theorem now follows at once from Proposition 4.3, making use of the fact that

$$\#\text{coker}(\text{Res}')/\#\ker(\text{Res}') = \#\Sha(A_{/K})^G/\#\Sha(A_{/F}). \quad \square$$

In the following corollary, we write $_N\Sha(A_{/K})$ for the kernel of the norm map $N_{K/F} : \Sha(A_{/K}) \to \Sha(A_{/K})^G$.

**Corollary 4.5.** *If conditions (A) and (B) above hold, then*

$$\#\hat{H}^0(G, \text{Ш}(A_{/K})) \cdot \#\text{Ш}(A_{/K}) = \#\,_N\text{Ш}(A_{/K}) \cdot \#\text{Ш}(A_{/F}) \cdot \prod_{v \in S} \#H^1(G_w, A(K_w)).$$

*If furthermore $K/F$ is a cyclic extension, then*

$$\#\text{Ш}(A_{/K}) = \#\,_N\text{Ш}(A_{/K}) \cdot \#\text{Ш}(A_{/F}).$$

*Proof.* The first assertion follows at once from the theorem and the exactness of the sequence

$$0 \to\,_N\text{Ш}(A_{/K}) \to \text{Ш}(A_{/K}) \xrightarrow{N_{K/F}} \text{Ш}(A_{/K})^G \to \hat{H}^0(G, \text{Ш}(A_{/K})) \to 0.$$

The second assertion follows from the first and the theorem, making use of the facts that, when $G$ is cyclic, $\#\hat{H}^0(G, \text{Ш}(A_{/K})) = \#H^1(G, \text{Ш}(A_{/K}))$ by [1], p. 109, and $\#H^1(G_w, A(K_w)) = \#H^2(G_w, A(K_w))$ if $w$ is non-archimedean by [15], §4 (14). □

The final considerations of this paper pertain to the case of quadratic extensions $K/F$, and are as follows.

Suppose that $K/F$ is a quadratic extension and let $\chi$ denote the non-trivial character of $G = \text{Gal}(K/F)$. We will write $A^\chi$ for the twist of $A$ by $\chi$ (see §2 of [9]). Then there is an isomorphism $\psi : A_{/K} \xrightarrow{\sim} A^\chi_{/K}$ such that $\psi^\sigma = \chi(\sigma)\psi$ for $\sigma \in G$. It follows that

(2) $$_N\text{Ш}(A^\chi_{/K}) \simeq \text{Ш}(A_{/K})^G.$$

**Corollary 4.6.** *Suppose that $K/F$ is a quadratic extension and let $\chi$ denote the non-trivial character of $G = \text{Gal}(K/F)$. Assume that conditions (A) and (B) above hold for both $A$ and $A^\chi$. Then*

$$\#\text{Ш}(A_{/K}) = \#\text{Ш}(A_{/F}) \cdot \#\text{Ш}(A^\chi_{/F}) \cdot \prod_{v \in S} \#H^1(G_w, A(K_w)).$$

Proof. By Corollary 4.5 applied to $A^\chi$ and (2), we have

$$\#\text{Ш}(A_{/K}) = \#\text{Ш}(A^\chi_{/K}) = \#\,_N\text{Ш}(A^\chi_{/K}) \cdot \#\text{Ш}(A^\chi_{/F})$$
$$= \#\text{Ш}(A_{/K})^G \cdot \#\text{Ш}(A^\chi_{/F}).$$

The corollary is now immediate from Theorem 4.4. □

Facultad de Ciencias, Universidad de Chile, Casilla 653, Santiago, Chile
*E-mail address*: cgonzale@abello.dic.uchile.cl